# FINITARY ISOMORPHISMS OF POISSON POINT PROCESSES

BY TERRY SOO[1] AND AMANDA WILKENS

*University of Kansas*

As part of a general theory for the isomorphism problem for actions of amenable groups, Ornstein and Weiss (*J. Anal. Math.* **48** (1987) 1–141) proved that any two Poisson point processes are isomorphic as measure-preserving actions. We give an elementary construction of an isomorphism between Poisson point processes that is finitary.

**1. Introduction.** We begin with some definitions and background necessary to state our main theorem. Let $r > 0$. A random variable $N$ taking values on $\mathbb{N} := \{0, 1, \ldots\}$ with $\mathbb{P}(N = m) = e^{-r} r^m / m!$ is a *Poisson* random variable with mean $r$. A *Poisson point process* on $\mathbb{R}^d$ with *intensity* $r$ is a random process $X$ taking values on the space $\mathbb{M}$ of Borel simple point measures (measures which are a countable sum of mutually singular delta measures) such that for every Borel subset $A \in \mathcal{B} := \mathcal{B}(\mathbb{R}^d)$ with finite Lebesgue measure $\mathcal{L}(A)$, the number of points of $X$ in $A$, denoted by $X(A)$, is a Poisson random variable with mean $r\mathcal{L}(A)$, and for any finite number of pairwise disjoint Borel sets $A_1, \ldots, A_\ell$ the random variables $X(A_1), \ldots, X(A_\ell)$ are independent.

Let $G$ be the group of isometries of $\mathbb{R}^d$. For each $r > 0$, let $P_r$ be the law of a Poisson point process with intensity $r$. We refer to the measure-preserving system $(\mathbb{M}, P_r, G)$, where $G$ acts on $\mathbb{M}$ via $g(\mu) = \mu \cdot g^{-1}$ for $g \in G$ and $\mu \in \mathbb{M}$, as a *Poisson system* with intensity $r$. A map $\phi : \mathbb{M} \to \mathbb{M}$ is a *factor* from $r$ to $s$ if on a set of $P_r$ full-measure, $\phi \circ g = g \circ \phi$ for all $g \in G$ and $P_r \circ \phi^{-1} = P_s$. A factor $\phi$ is an *isomorphism* if it is a bijection almost surely, in which case $\phi^{-1}$ serves as a factor from $s$ to $r$.

Ornstein and Weiss proved the following theorem in [12], as part of a more general theory.

THEOREM 1 (Ornstein and Weiss). *Any two Poisson systems are isomorphic.*

In particular, Theorem 1 follows from [12], Theorem 2, page 117, and verifying that Poisson systems are *extremal*, a technical condition; Ornstein and Weiss verified that systems associated with Poisson point processes on $\mathbb{R}^d$ endowed with the (smaller) group of *translations* of $\mathbb{R}^d$ are extremal. We also note that in the

Received May 2018; revised December 2018.
[1]Supported in part by New Faculty General Research Fund and General Research Fund.
*MSC2010 subject classifications.* 37A35, 60G10, 60G55.

special case where $d = 1$, a Poisson system is a canonical example of an infinite entropy Bernoulli flow and a translation-equivariant version of Theorem 1 is an consequence of Ornstein theory [11]. See also [14] for background and related problems for the case of Bernoulli flows.

In our paper, we give a proof of Theorem 1 by constructing an explicit isomorphism in the proof of our main theorem (Theorem 2), and we gain a nice property for the isomorphism map in the process.

We use the notation $\mu|_A(\cdot) := \mu(\cdot \cap A)$ for the restriction of $\mu$ to $A \in \mathcal{B}$, $B(z, \varepsilon) \subseteq \mathbb{R}^d$ for the open Euclidean ball of radius $\varepsilon$ centered at $z \in \mathbb{R}^d$, and $\mathbf{0}$ for the origin in $\mathbb{R}^d$. For $z \in \mathbb{R}^d$, we let $t_z \in G$ denote translation by $z$. Let $\phi$ be a factor map from $r$ to $s$, and let $\mu, \mu' \in \mathbb{M}$. We say a *coding window* of $\phi$ is a function $w : \mathbb{M} \to \mathbb{N} \cup \{\infty\}$ such that if $\mu|_{B(\mathbf{0}, w(\mu))} = \mu'|_{B(\mathbf{0}, w(\mu))}$, then $\phi(\mu)|_{B(\mathbf{0}, 1)} = \phi(\mu')|_{B(\mathbf{0}, 1)}$.

We say that $\phi$ is *finitary* if there exists a coding window $w$ such that $w$ is finite $P_r$-almost surely. Since $\phi$ is translation-equivariant, we have $\phi(\mu)|_{B(z, 1)} = \phi(t_z^{-1}\mu)|_{B(\mathbf{0}, 1)}$ for any $z \in \mathbb{R}^d$, so that if $\phi$ is finitary, then the values of $\phi(\mu)$ restricted to any unit ball are determined by the values of $\mu$ restricted to a (larger) concentric ball.

THEOREM 2 (Finitary isomorphism). *Any two Poisson systems are isomorphic; furthermore, there exists a finitary isomorphism with finitary inverse.*

The problem of determining when two measure-preserving systems are isomorphic has a long history [3, 16] and the most understood systems are those associated with independent and identically distributed (i.i.d.) processes indexed by the integers, where Kolmogorov–Sinai entropy is a complete isomorphism invariant. Ornstein proved that any two equal entropy i.i.d. processes are isomorphic as factors [10] and Keane and Smorodinsky strengthened this result by constructing almost everywhere continuous isomorphisms between any two processes of finite equal entropy [8, 9]. In this discrete setting, the continuity of the factor map is equivalent to the factor map possessing a random finite coding window, which is what we adapt to be the definition of finitary in the point process setting; in the point process setting, it may appear to be a strong requirement, since one can determine the *exact* location of all the points of the output process in the unit ball, given a large enough coding window.

We say that $U$ is a *uniform random variable* if it is uniformly distributed on the unit interval $[0, 1]$. Kalikow and Weiss proved that when the group of isometries of a Poisson system on the real line is restricted to translations by a unit length, then the Poisson system is finitarily isomorphic to the infinite entropy Bernoulli shift given by independent uniform random variables indexed by the integers [7]; it is then immediate that, in this restricted case, any two such Poisson systems on the real line are finitarily isomorphic.

Our proof of Theorem 2 will make use of a key construction due to Holroyd, Lyons, and Soo in [6], wherein they proved any two Poisson systems are finitarily homomorphic in the following sense.

THEOREM 3 (Holroyd, Lyons and Soo). *Fix $s > 0$. There exists $\phi : \mathbb{M} \to \mathbb{M}$ such that for all $r > 0$, the map $\phi$ is a finitary factor from $r$ to $s$.*

As in [6], when we build a map to generate a Poisson point process from a Poisson point process, we use randomness harnessed from the input system in a careful way so as not to disrupt independence of the system. Once independence is assured, we convert the randomness to a uniform random variable, and then convert the uniform random variable to a Poisson point process on a finite volume (specifically, a cell of an isometry-equivariant partition). At each step our maps will be entirely explicit. We remark that an injective measurable map and thus isomorphism from a uniform random variable to any Poisson point process on a finite volume cannot exist since the unique empty point measure, which we denote by $\varnothing$, occurs with nonzero probability.

To circumvent the nonexistence of such an isomorphism, we prove Proposition 4 as an intermediate result, from which Theorem 2 will follow. To state Proposition 4, we need a few more definitions.

Let $X$ and $Y$ be independent Poisson point processes on $\mathbb{R}^d$ with respective intensities $r > 0$ and $s > 0$ and let $\psi : \mathbb{M} \to \mathbb{M} \times \mathbb{M}$. For $\mu \in \mathbb{M}$ we write $\psi(\mu) = (\psi(\mu)_1, \psi(\mu)_2)$. We say that $\psi$ is *isometry-equivariant* if

$$\psi(g\mu) = g\psi(\mu) := (g\psi(\mu)_1, g\psi(\mu)_2)$$

for all $\mu \in \mathbb{M}$ and all isometries $g$ of $\mathbb{R}^d$. If $\psi(X)_1$ is a Poisson point process on $\mathbb{R}^d$ of intensity $r$ independent of $\psi_s(X)$, a Poisson point process on $\mathbb{R}^d$ of intensity $s$, and on a set of $P_r$ full-measure, the map $\psi$ is isometry-equivariant, then we say that $\psi$ is a *factor from $r$ to $(r, s)$*. The map $\psi_s$ is *finitary* if each coordinate mapping is finitary.

Again, let $r, s > 0$, and let $X$ and $Y$ be as above. Let $\zeta : \mathbb{M} \times \mathbb{M} \to \mathbb{M}$. We say $\zeta$ is finitary if there exists a coding window $w : \mathbb{M} \times \mathbb{M} \to \mathbb{N} \cup \{\infty\}$ such that $w$ is finite $(P_r \times P_s)$-almost surely and we have that if

$$(\mu_1, \mu_2)|_{B(\mathbf{0}, w(\mu_1, \mu_2))} = (\mu'_1, \mu'_2)|_{B(\mathbf{0}, w(\mu_1, \mu_2))},$$
$$\text{then } \zeta(\mu_1, \mu_2)|_{B(\mathbf{0}, 1)} = \zeta(\mu'_1, \mu'_2)|_{B(\mathbf{0}, 1)}.$$

If $\zeta(X, Y)$ is a Poisson point process on $\mathbb{R}^d$ of intensity $r$, and the map $\zeta$ is isometry-equivariant, then we say that $\zeta$ is a *factor from $(r, s)$ to $r$*. Thus a factor $\psi$ from $r$ to $(r, s)$ is an *isomorphism* if it is a bijection almost surely, in which case its inverse serves as a factor from $(r, s)$ to $r$.

PROPOSITION 4. *Fix $s > 0$. There exists a finitary isomorphism $\psi_s$ from $r$ to $(r, s)$, independent of $r > 0$. Furthermore, the map $\psi_s$ has finitary inverse.*

The map $\psi_s$ applied to a Poisson point process of intensity $r$ yields two Poisson point processes, one of intensity $r$ and one of $s$. The process of intensity $r$ differs from the original process only within particular unit balls, each of which contains a unique point of the original process—using the randomness of these points we resample them and generate the Poisson point process of intensity $s$. The additional information contained within the process of intensity $r$ allows $\psi_s$ to be injective. Care is required to ensure that we do not violate independence properties within each system.

Of course, after an application of $\psi_s$ we are left with too much information rather than not enough for an isomorphism between two singular Poisson systems. We make slight adjustments in the proof of Theorem 2, in Section 3.3. As an additional, rather immediate consequence of Proposition 4, we have that Poisson systems are finitarily isomorphic to products of Poisson systems (Theorem 25).

**2. Preliminary results.** We work toward a constructive proof of Theorem 2, utilizing the framework found in [6]. We will refer to the restriction of a Poisson point process on $\mathbb{R}^d$ to $A \in \mathcal{B}$ as a *Poisson point process on A*. One key idea is to use some randomness of the input Poisson point process to obtain an isometry-equivariant partition of $\mathbb{R}^d$ and then to generate Poisson point processes of some desired intensity on each cell of the partition, thus yielding a Poisson point process of the same intensity on the whole of $\mathbb{R}^d$. That we indeed end up with a such a process on $\mathbb{R}^d$ is immediate from Remark 5 below, as long as we are careful to satisfy required independence properties.

2.1. *Uniform random variables.* Recall $U$ is a uniform random variable if it is uniformly distributed on the unit interval. Similarly, we say that $U$ is a *uniform random variable on A* if $\mathbb{P}(U \in \cdot) = \mathcal{L}(\cdot \cap A)/\mathcal{L}(A)$. Often we use the notation $U[A]$ to denote a uniform random variable on $A$.

REMARK 5 (Uniformly distributed random variables and Poisson point processes). Throughout the paper we use that if $X$ is a Poisson point process on $\mathbb{R}^d$ and $A \in \mathcal{B}$ nonempty with finite Lebesgue measure, then conditional on the event that $X(A) = n$, these $n$ points of $X$ are independently and uniformly distributed in $A$; this property in fact characterizes Poisson point processes.

Thus if $\{U[A]_i\}_{i \in \mathbb{N}}$ is an i.i.d. sequence of uniform random variables on $A$ and $N$ is a Poisson random variable with mean $r\mathcal{L}(A)$ independent of the sequence, we may write

$$(2.1) \qquad X|_A \stackrel{d}{=} \sum_{i=1}^{N} \delta[U[A]_i],$$

where $\stackrel{d}{=}$ denotes equality in distribution, and we denote Dirac measure with mass at $z$ as $\delta[z] \in \mathbb{M}$.

Let $A_1, A_2, \ldots \in \mathcal{B}$ with finite Lebesgue measure such that $\{A_i\}_{i\in\mathbb{N}}$ partitions $\mathbb{R}^d$. From elementary properties of Poisson point processes, the random variables $X|_{A_1}, X|_{A_2}, \ldots$ are independent. Thus generating independent Poisson point processes on the cells of the partition via (2.1) generates a Poisson point process on $\mathbb{R}^d$.

In Example 6 below, we illustrate one way to generate a Poisson point process on $\mathbb{R}^2$ given a partition and a uniform random variable for each cell.

EXAMPLE 6 (Generating a Poisson point process on $\mathbb{R}^2$). Let $(k_1, k_2) = k \in \mathbb{Z}^2$, and let $C_k$ be the unit square with set of endpoints

$$\{(0,0)+k, (1,0)+k, (1,1)+k, (0,1)+k\}$$

so that $\{C_k\}_{k\in\mathbb{Z}^2}$ is almost surely a partition of $\mathbb{R}^2$. Let $U = \{U_k\}_{k\in\mathbb{Z}^2}$ be a collection of independent uniform random variables. We construct a family of measurable maps $\{\pi_k^r\}_{k\in\mathbb{Z}^2}$ where each map $\pi_k^r : [0,1] \to \mathbb{M}$ sends $U_k$ to a Poisson point process with intensity $r > 0$ on $C_k$. To do so each map will perform several actions; namely, determining the number of points of $\pi_k^r(U_k)$ in $C_k$ and generating each point in such a way so that $\pi_k^r(U_k)$ has distribution as in (2.1). Set $\pi_k^r = \pi_k$.

Let $N$ be a Poisson random variable with mean $r$ and set $p_m = \mathbb{P}(N < m)$. Note $\{[p_m, p_{m+1})\}_{n=0}^\infty$ partitions $[0,1)$. Define $q : [0,1) \to \mathbb{N}$ piecewise so that whenever $x \in [p_m, p_{m+1})$ then $q(x) = m$ for $0 \le m < \infty$. Also define $f : [0,1) \to [0,1)$ piecewise so that whenever $x \in [p_m, p_{m+1})$,

$$f(x) = \frac{x - p_m}{p_{m+1} - p_m}.$$

Conditioned on $q(U_k) = m$, we have that $U_k$ is a uniform random variable on $[p_m, p_{m+1})$, hence

$$\mathbb{P}\big(f(U_k) \in A | q(U_k) = m\big) = \mathbb{P}(U_k \in A)$$

for any $A \in \mathcal{B}([0,1])$ and any $m \in \mathbb{N}$. Thus $f(U_k) \stackrel{d}{=} U_k$ and is independent of $q(U_k)$. We use the randomness of $f(U_k)$ to populate points in $C_k$. To ensure each point's location is independent of all other points' locations, we split $f(U_k)$ into distinct independent uniform random variables.

For $x \in [0,1]$, let $.x_1 x_2 x_3 \ldots$ be the binary expansion of $x$. Define $b_n : [0,1] \to [0,1]^n$ so that for $x \in [0,1]$ we have

(2.2) $$b_n(x) = (x^1, x^2, \ldots, x^n),$$

where $x^i = .x_i x_{n+i} x_{2n+i} \ldots$ for $1 \le i \le n$. We denote the $i$th coordinate of $b_n(x)$ as $b_n(x)^i$. We will apply $b_n$ for some $n$ to $f(U_k)$. Set

$$\pi_k(U_k) = \sum_{i=1}^{q(U_k)} \delta\big[\big(b_{2q(U_k)}(f(U_k))^i + k_1, b_{2q(U_k)}(f(U_k))^{2q(U_k)-i+1} + k_2\big)\big]$$

if $q(U_k) \neq 0$. Otherwise, set $\pi_k(U_k) = \varnothing$. Each pair

$$\left(b_{2q(U_k)}(f(U_k))^i + k_1, b_{2q(U_k)}(f(U_k))^{2q(U_k)-i+1} + k_2\right)$$

is a uniformly distributed point in $C_k$.

Define $\pi : [0,1]^{\mathbb{Z}^2} \to \mathbb{M}$ so that $\pi(U)|_{C_k} = \pi_k(U_k)$ for all $k \in \mathbb{Z}$. Then $\pi(U)$ is a Poisson point process of intensity $r$ on $\mathbb{R}^2$ by Remark 5.

Note the mapping in Example 6 is not isometry-equivariant and is merely an indication of how one might generate a Poisson point process on $\mathbb{R}^2$. To simplify our approach to proving Proposition 4, we will first prove the translation-equivariant version (Proposition 15).

Another key idea arises from the problem of injectivity. In Example 6, we were careful not to let any information go to waste. We could have easily generated up to infinitely many uniform random variables from the first by a function similar to $b_n$ in (2.2) and assigned one to provide the number of points of $\pi_k(U_k)$ in $C_k$. The map in our example is closer to being injective than such a map would be, but it is easy to see where injectivity fails—there are infinitely many ways to obtain the outcome $q(U_k) = 0$ but only one empty point process.

2.2. *An isometry-equivariant partition.* We construct our desired isometry-equivariant partition of $\mathbb{R}^d$ in two phases. The general idea for our first phase is to partition $\mathbb{R}^d$ into two sets; the one containing balls of a certain type, and the other containing everything else. The following definitions match those in [6].

Let $X$ be a Poisson point process on $\mathbb{R}^d$ with intensity $r > 0$. Define a *shell* centered at $x$ from $a$ to $b$ as the set

$$L(x,a,b) = \{y \in \mathbb{R}^d : a \leq \|x-y\| \leq b\}.$$

Recall that $X(A)$ is the number of points of $X$ in $A$. We call a point $x \in \mathbb{R}^d$ a *pre-seed* if $X(L(x, 78+d, 105+d)) = 0$ and for every open ball $B$ of radius 0.5 contained strictly inside $L(x, 11, 78+d)$, we have $X(B) \geq 1$. Although the probability that $B(\mathbf{0}, 1)$ contains a pre-seed is small, infinitely many $x \in \mathbb{R}^d$ are pre-seeds $P_r$-almost surely. If $x$ is a pre-seed, then we refer to $L(x, 78+d, 105+d)$ as its *empty shell* and $L(x, 11, 78+d)$ as its *halo*. Figure 1 (which also appears in [6]) illustrates a pre-seed.

Given two pre-seeds $x$ and $y$, by definition either $\|x-y\| \leq 2$ or $\|x-y\| \geq 132+d$, since the empty shell of $x$ cannot intersect the halo of $y$ and vice versa. We say that $x$ and $y$ are *related* if $\|x-y\| \leq 2$. Thus we have an equivalence relation on pre-seeds. Let $C$ be an equivalence class of pre-seeds under $X$, so that $C$ is contained in some ball of radius 2. Then there exists a ball containing $C$ with unique smallest radius and center $c$. We say that $c$ is a *seed*. Although $c$ may not be a pre-seed, we still refer to $L(c, 78+d, 105+d)$ as its halo. Using seeds we can precisely define our first partition.

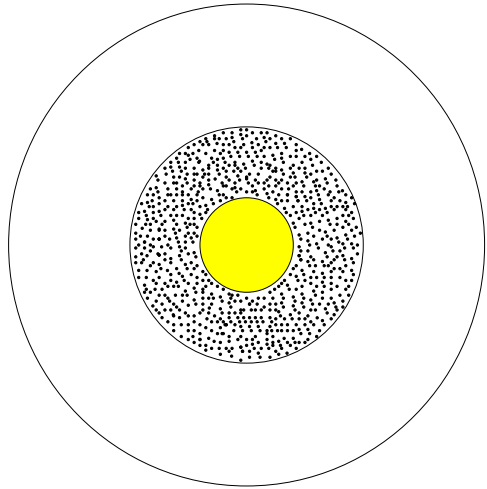

FIG. 1. *A pre-seed. The empty shell contains no points of X. The halo is relatively densely filled with points of X. The shaded area is unspecified in terms of X.*

For every seed $c$, we call $\bar{B}(c, 1)$ a *globe*. Let $\mathcal{F}$ be the set of closed subsets of $\mathbb{R}^d$. Define $\mathcal{S} : \mathbb{M} \to \mathcal{F}$ so that for $\mu \in \mathbb{M}$, we have that $\mathcal{S}(\mu)$ is the union of the set of globes under $\mu$. We say that $\mu, \mu' \in \mathbb{M}$ *agree* on a set $A \in \mathcal{B}$ if their restrictions to $A$ are equal. By [6], Proposition 15 and Lemma 32, the mapping $\mathcal{S}$ has the following properties:

(a) If $X$ is a Poisson point process on $\mathbb{R}^d$ with intensity $r$ and law $P_r$, then $P_r$-almost surely $\mathcal{S}(X)$ is a nonempty union of disjoint closed balls of radius 1.

(b) The map $\mathcal{S}$ is isometry-equivariant; that is, for all $g \in G$ and $\mu \in \mathbb{M}$, $\mathcal{S}(g\mu) = g\mathcal{S}(\mu)$.

(c) For all $\mu, \mu' \in \mathbb{M}$, if $\mu$ and $\mu'$ agree on the set

$$\left( \bigcup_{x \in \mathcal{S}(\mu)} \bar{B}(x, 2) \right)^c,$$

then $\mathcal{S}(\mu) = \mathcal{S}(\mu')$.

(d) Furthermore, for any $z \in \mathbb{R}^d$ any $\mu, \mu' \in \mathbb{M}$, if $\bar{B}(z, 1)$ is a globe under $\mu$, then whenever $\mu$ and $\mu'$ agree on $B(z, 125 + d)$, then $\bar{B}(z, 1)$ is also a globe under $\mu'$.

We refer to $\mathcal{S}$ as a *selection rule*. We denote the set of globes by Globes$[\mathcal{S}(X)]$. Note Globes$[\mathcal{S}(X)]$ and $\mathcal{S}(X)^c$ partition $\mathbb{R}^d$ in an isometry-equivariant way. This first phase of our partition allows us to harness the randomness that we need in order to generate a Poisson point process from a given Poisson point process. Consider again property (c). That $\mathcal{S}(X)$ does not hold information on the Poisson point process within the globes is an important distinction. It is also important to note $\mathcal{S}(X)$ depends only on $X$ restricted to $\mathcal{S}(X)^c$ (in fact, slightly less than $\mathcal{S}(X)^c$ by the definition; this is a relic of the proof of [6], Proposition 15).

Property (d) is a localized version of property (c) we will use to ensure the isomorphism we define for Theorem 2 is finitary.

If a globe contains a unique point of $X$, we call the globe *special*, and we let $\mathcal{S}^*(X)$ be the union of the special globes. We denote the set of special globes by Globes*$[\mathcal{S}(X)]$. (By the upcoming Proposition 8 there are infinitely many special globes $P_r$-almost surely.) Since points of $X$ are uniformly distributed inside any nonempty finite volume Borel subset, we may think of the point in a special globe as a uniform random variable on a closed ball of radius 1. In Lemma 7, we detail an explicit map from a closed ball of radius 1 to the unit interval.

We use the following facts in the proof of Lemma 7. A nonzero point in $\mathbb{R}^d$ may be written uniquely in polar coordinates as

$$(r, \theta_1, \ldots, \theta_{d-2}, \theta_{d-1}),$$

where $r$ is the distance to the origin, the angles $\theta_1, \ldots, \theta_{d-2}$ range from 0 to $\pi$, and the angle $\theta_{d-1}$ ranges from 0 to $2\pi$.

Let $(R, \Theta_1, \ldots, \Theta_{d-2}, \Theta_{d-1})$ be uniformly distributed in the closed ball $\bar{B}(\mathbf{0}, 1)$ and let $F_0, F_1, \ldots, F_{d-2}, F_{d-1}$ be the cumulative distribution functions (cdf's) for the random variables $R, \Theta_1, \ldots, \Theta_{d-2}, \Theta_{d-1}$, respectively. Note the cdf's are continuous and increasing, so their inverse functions are well defined. Thus if $U$ is a uniform random variable, then

$$F_0(R) \stackrel{d}{=} F_1(\Theta_1) \stackrel{d}{=} \cdots \stackrel{d}{=} F_{d-1}(\Theta_{d-1}) \stackrel{d}{=} U.$$

Indeed, one may check that $F_0(R) = R^d$ for any $d$, but it becomes difficult (in the sense of integrating powers of trigonometric functions) to write $F_i(\Theta_i)$ explicitly for high dimensions.

Furthermore, the random variables $R, \Theta_1, \ldots, \Theta_{d-2}, \Theta_{d-1}$ are independent.

LEMMA 7 (Ball to unit interval isomorphism). *For every $d \geq 1$ there exists an isomorphism $b^d : \bar{B}(\mathbf{0}, 1) \to [0, 1]$ such that if $V$ is uniformly distributed on $\bar{B}(\mathbf{0}, 1)$, then $b^d(V)$ is a uniform random variable on $[0, 1]$.*

PROOF. We prove the case $d = 1$ separately and first let $d \geq 2$. Let $z \in \bar{B}(\mathbf{0}, 1)$ and write the polar coordinates of $z$ as $(r, \theta_1, \ldots, \theta_{d-2}, \theta_{d-1})$. Suppose $(R, \Theta_1, \ldots, \Theta_{d-2}, \Theta_{d-1})$ is uniformly distributed in $\bar{B}(\mathbf{0}, 1)$ and let $F_0$ be the cumulative distribution function for $R$ and $F_i$ the cumulative distribution function for $\Theta_i$, where $1 \leq i \leq d - 1$.

Recall notation for the binary expansion of an element in $[0, 1]$ (as used in (2.2)) which we use in the following definition. Let $b^d : \bar{B}(\mathbf{0}, 1) \to [0, 1]$ such that

$$b^d(z) = .[F_0(r)]_1[F_1(\theta_1)]_1 \ldots [F_{d-2}(\theta_{d-2})]_1[F_{d-1}(\theta_{d-1})]_1[F_0(r)]_2 \ldots.$$

The map $b^d$ interweaves the binary expansions of the polar coordinates $F_0(r)$, $F_1(\theta_1), \ldots, F_{d-2}(\theta_{d-2}), F_{d-1}(\theta_{d-1})$ and outputs a single element in $[0, 1]$. Since the coordinates are independent, if $V$ is a $U[\bar{B}(\mathbf{0}, 1)]$ random variable then $b^d(V)$ is a uniform random variable on $[0,1]$. The map $b^d$ is bijective almost surely.

In the case $d = 1$, we have no need for polar coordinates. Let $z \in \bar{B}(\mathbf{0}, 1)$. Define $b^1 : \bar{B}(\mathbf{0}, 1) \to [0, 1]$ so that

$$b^1(z) = \frac{z+1}{2}.$$

It is simple to check $b^1$ satisfies our conditions. □

Lemma 7 provides us with a mechanism to extract uniform random variables from a Poisson point process via the special globes. A selection rule $\mathcal{S}$ implies that such uniform random variables are conditionally independent of the process outside of the globes. Now we make concrete our key idea of independence using Proposition 16 from [6], appearing here as Proposition 8.

PROPOSITION 8 (Holroyd, Lyons, and Soo). *Let $X$ and $W$ be independent Poisson point processes on $\mathbb{R}^d$ with the same intensity. For a selection rule $\mathcal{S}$, the process $Z := W|_{\mathcal{S}(X)} + X|_{\mathcal{S}(X)^c}$ has the same law as $X$ and $\mathcal{S}(X) = \mathcal{S}(Z)$.*

Thus for a Poisson point process $X$ and a selection rule $\mathcal{S}$, given the knowledge of $\mathcal{S}(X)$, we have that $X|_{\mathcal{S}(X)}$ is a Poisson point process on $\mathcal{S}(X)$ independent of $X|_{\mathcal{S}(X)^c}$. In particular, we rely on the fact that knowing $X|_{\mathcal{S}(X)^c}$ does not give us any information on the location of points inside the globes.

So that we may reference elements in Globes*$[\mathcal{S}(X)]$, we let $\{\beta_i\}_{i\in\mathbb{N}} = $ Globes*$[\mathcal{S}(X)]$ where the special globes are ordered by the distance of their centers to the origin. Let $c_i$ be the center of the special globe $\beta_i$ and $x_i$ the unique point of $X$ in $\beta_i$. We list two applications of Proposition 8 that we will use in the proof of Proposition 4.

COROLLARY 9. *Let $X$ be a Poisson point process on $\mathbb{R}^d$ with intensity $r > 0$. Let $\{U_i[B]\}_{i\in\mathbb{N}}$ be a sequence of independent uniform random variables on $B = \bar{B}(\mathbf{0}, 1)$ that is independent of $X$. Then*

$$\left(X|_{\mathcal{S}^*(X)^c}, \mathcal{S}^*(X), X|_{\mathcal{S}^*(X)}\right) \stackrel{d}{=} \left(X|_{\mathcal{S}^*(X)^c}, \mathcal{S}^*(X), \sum_{i\in\mathbb{N}} \delta[t_{c_i}(U_i[B])]\right).$$

PROOF. Immediate from Proposition 8 and Remark 5. □

COROLLARY 10. *Let $X$ be a Poisson point process on $\mathbb{R}^d$ with intensity $r > 0$. There exists a measurable map $\alpha : \mathbb{M} \times \mathcal{F} \to [0, 1]$ such that if $\{U_i\}_{i\in\mathbb{N}}$ is a sequence of independent uniform random variables that is independent of $X$, then*

$$\left(X|_{\mathcal{S}^*(X)^c}, \mathcal{S}^*(X), \{\alpha(X, \beta_i)\}_{i\in\mathbb{N}}\right) \stackrel{d}{=} \left(X|_{\mathcal{S}^*(X)^c}, \mathcal{S}^*(X), \{U_i\}_{i\in\mathbb{N}}\right).$$

PROOF. Let $b = b^d$ be the map from Lemma 7. Define $\alpha : \mathbb{M} \times \mathcal{F} \to [0, 1]$ so that for $\mu \in \mathbb{M}$, $\beta \in \mathcal{F}$, and $\beta_i \in \text{Globes}^*[\mathcal{S}(\mu)]$, we have

$$\alpha(\mu, \beta) = \begin{cases} b(x_i - c_i) & \text{whenever } \beta = \beta_i, \\ 0 & \text{otherwise.} \end{cases}$$

By Corollary 9, we have

$$(X|_{\mathcal{S}^*(X)^c}, \mathcal{S}^*(X), \{t_{c_i}^{-1}(x_i)\}_{i \in \mathbb{N}}) \stackrel{d}{=} (X|_{\mathcal{S}^*(X)^c}, \mathcal{S}^*(X), \{U_i[B]\}_{i \in \mathbb{N}}).$$

Let $\{U_i\}_{i \in \mathbb{N}}$ be a sequence of independent uniform random variables independent of $X$. Note $t_{c_i}^{-1}(x_i)$ takes values in $\bar{B}(\mathbf{0}, 1)$. Thus, we have

$$(X|_{\mathcal{S}^*(X)^c}, \mathcal{S}^*(X), \{\alpha(X, \beta_i)\}_{i \in \mathbb{N}}) \stackrel{d}{=} (X|_{\mathcal{S}^*(X)^c}, \mathcal{S}^*(X), \{U_i\}_{i \in \mathbb{N}}),$$

after an application of $b$. □

Given a Poisson point process, we now have a way to extract randomness within the process carefully enough to respect independence. The isomorphism we are working to construct will make good use of this randomness, but first we need the second phase of our isometry-equivariant partition.

Let $X$ be a Poisson point process on $\mathbb{R}^d$ with intensity $r$, special globes $\beta_i$, and centers of the special globes $c_i$. Now we refer to each $c_i$ as a *site*. The *special Voronoi cell* of a site $c_i$ is the set of all points $y \in \mathbb{R}^d$ such that $\|c_i - y\| < \|y - c_k\|$ for all $i \neq k$. Remark 11 follows from our definitions and will be used in the proof of the finitary property.

REMARK 11. The law of the point process of sites given by $\sum_{i \in \mathbb{N}} \delta[c_i]$ is translation-invariant, and thus the special Voronoi cells are bounded convex polytopes. In addition, the Voronoi cell of a site $c$ contains the globe $\bar{B}(c, 1)$ and its halo.

We define the *special Voronoi tessellation* $\mathcal{V}^*(X)$ to be the set of special Voronoi cells of sites $c_i$ for all $i \in \mathbb{N}$. Note $\mathcal{V}^*(X)$ partitions $\mathbb{R}^d$ $P_r$-almost surely. It is clear that $\mathcal{V}^*(X)$ is itself isometry-equivariant; for any isometry $g \in G$, we have

$$\mathcal{V}^*(gX) = \{gv : v \in \mathcal{V}^*(X)\} = g\mathcal{V}^*(X).$$

Our isomorphism will output Poisson point processes of desired intensity via uniform random variables gathered from the input Poisson point process, within each cell of the special Voronoi tessellation.

**3. Proof of Theorem 2.** We have introduced our isometry-equivariant partition and methods for extracting randomness. We need further tools to establish how we will obtain a finitary, injective map between Poisson point processes.

3.1. *Tools for the finitary property and injectivity.* We prove that the special Voronoi cells of a Poisson point process are locally determined in Lemma 12, adapted from [6], Theorem 31. For $\mu \in \mathbb{M}$, we let $v(\mu, z)$ be the Voronoi cell such that $z \in v(\mu, z)$.

LEMMA 12 (Local property of special Voronoi cells). *Let $r > 0$. There exists a map $w : \mathbb{M} \to \mathbb{N} \cup \{\infty\}$ such that $w$ is finite $P_r$-almost surely and for $P_r$-almost all $\mu, \mu' \in \mathbb{M}$, if $\mu$ and $\mu'$ agree on $B(\mathbf{0}, w(\mu))$, then for all $z \in B(\mathbf{0}, 1)$ we have $v(\mu, z) \subseteq B(\mathbf{0}, w(\mu))$ and $v(\mu, z) = v(\mu', z)$.*

PROOF. Let $\mu, \mu' \in \mathbb{M}$. Recall that the sites are centers of special globes, and by property (d) the globes are locally determined in the following sense: if $\bar{B}(z, 1)$ is a special globe under $\mu$ and $\mu$ agrees with $\mu'$ on a sufficiently large ball about $z$, then $\bar{B}(z, 1)$ is also a special globe under $\mu'$. Thus, it suffices to find the radius of a ball containing sufficiently many sites to determine the Voronoi cells that intersect $\bar{B}(\mathbf{0}, 1)$.

Set $\ell = 100(106 + d)$. Let $\{C_k\}_{k \in \mathbb{Z}^d}$ partition $\mathbb{R}^d$ into equal sized cubes of side length $\ell$ so that $C_k$ is centered at $k\ell$. Then $B(k\ell, 1) \subset C_k$. Let $E_k$ be the event that $B(k\ell, 1)$ contains the center of a special globe. By Proposition 8, the $E_k$ are independent under the Poisson measure $P_r$ and occur with nonzero probability.

Let $T_1(\mu)$ be the smallest integer $n$ such that there exists integers $k_i$ such that

$$-n < k_{-3} < k_{-2} < k_{-1} < 0 < k_1 < k_2 < k_3 < n$$

and events $E_{(k_i, 0, \ldots, 0)}$ all occur. For each coordinate $i = 1, \ldots, d$, we similarly define $T_i$. Now set $w = 8\ell \sum_{i=1}^{d} T_i$. Any Voronoi cell intersecting $B(\mathbf{0}, 1)$ is contained in $B(\mathbf{0}, w(\mu))$ for $P_r$ almost all $\mu$, and all such Voronoi cells are determined by restriction of $\mu$ to this ball. Moreover, it is easy to verify that if $X$ is a Poisson point process of intensity $r$, then $\mathbb{E}w(X) < \infty$. □

Let $\mathcal{K} \subset \mathcal{B}(\mathbb{R}^d)$ denote the set of bounded convex polytopes of dimension $d$. Recall that by Remark 11, the Voronoi cells of the special Voronoi tessellation of a Poisson point process with law $P_r$ are $P_r$-almost surely elements of $\mathcal{K}$.

Also recall that there exists a measurable map from a uniform random variable to any Poisson point process on a finite volume, but an isomorphism cannot exist since the unique point measure $\varnothing$ occurs with nonzero probability. To circumvent this obstruction to our isomorphism, in Proposition 14, given a uniform random variable we generate a new uniform random variable in addition to a Poisson point process for each element of $\mathcal{K}$.

First, we state Lemma 13. We will use the map given here to construct a Poisson point process in an arbitrary element of $\mathcal{K}$.

LEMMA 13. *Let $n \geq 1$, and let $A \in \mathcal{K}$. There exists a measurable injection $g : [0, 1] \to \mathbb{M}$ such that if $U$ is a uniform random variable, then*

$$g(U) \stackrel{d}{=} \sum_{i=1}^{n} \delta[U_i[A]],$$

*where $U_1[A], \ldots, U_n[A]$ are independent and uniformly distributed on $A$.*

We give a constructive and elementary proof of Lemma 13 in Section 4.1; one may also refer to a version of the Borel isomorphism theorem for standard probability spaces, at the cost of concreteness. See [15], Theorem 3.4.23.

PROPOSITION 14. *Let $r > 0$. There exists a collection of measurable maps $\{\pi_{(A,r)}\}_{A \in \mathcal{K}}$ where for each $A \in \mathcal{K}$, the map $\pi_A := \pi_{(A,r)} : [0, 1] \to \mathbb{M} \times [0, 1]$ has the following properties. Write $\pi_A(U) = (\pi_A(U)_1, \pi_A(U)_2)$.*

(a) *If $U$ is a uniform random variable, then $\pi_A(U)_1$ is a Poisson point process of intensity $r$ on $A$ and $\pi_A(U)_2$ is a uniform random variable.*
(b) *The Poisson point process $\pi_A(U)_1$ is independent of the uniform random variable $\pi_A(U)_2$.*
(c) *Each map $\pi_A$ is injective almost surely.*

PROOF. Fix $A \in \mathcal{K}$. Let $U$ be a uniform random variable. By Lemma 13, for each $m \geq 1$, let $g_m : [0, 1] \to \mathbb{M}$ be a measurable injection so that $g_m(U)$ has the law of $m$ independent random variables uniformly distributed on $A$. We define some functions similar to those in Example 6. Let $N$ be a Poisson random variable with mean $r\mathcal{L}(A)$. For each $m \in \mathbb{N}$, let $p_m = \mathbb{P}(N < m)$. Note that $p_0 = 0$. Let $x \in [0, 1)$. Define $q$ and $f$ so that for $x \in [p_m, p_{m+1})$ we have

$$q(x) = m \quad \text{and} \quad f(x) = \frac{x - p_m}{p_{m+1} - p_m}.$$

Then $f(U)$ is a uniform random variable independent of $q(U)$ as in Example 6.

Let $b_2$ be the binary expansion map in (2.2) with $n = 2$, so we have $b_2 : [0, 1] \to [0, 1] \times [0, 1]$ where $b_2(x) = (x^1, x^2)$; we denote the $i$th coordinate of $b_2(x)$ as $b_2(x)^i$. In particular, we have that $b_2(U)^1$ and $b_2(U)^2$ are independent uniform random variables, and $b_2$ is injective almost surely.

Define $\pi_A : [0, 1] \to \mathbb{M} \times [0, 1]$ so that

$$\pi_A(x) = (\pi_A(x)_1, \pi_A(x)_2)$$

$$= \begin{cases} (\varnothing, f(x)) & \text{whenever } q(x) = 0, \\ (g_{q(x)}[b_2(f(x))^1], b_2(f(x))^2) & \text{otherwise.} \end{cases}$$

By Remark 5, $\pi_A(U)_1$ is a Poisson point process of intensity $r$ on $A$; moreover, $\pi_A(U)_1$ is independent of $\pi_A(U)_2$, a uniform random variable.

As for injectivity, in the case $\pi_A(x)_1 = \emptyset$, we have $q(x) = 0$. Given this along with $f(x)$, which we have as $\pi_A(x)_2$, we may reconstruct $x$ precisely. If $\pi_A(x)_1$ contains $m \geq 1$ points, we have $q(x) = m$. Since $g_m$ is injective, we recover $b_2(f(x))^1$. From $b_2(f(x))^1$, $b_2(f(x))^2$, and $m$ we recover $x$. Thus $\pi_A$ is injective almost surely. □

We are nearly ready to prove the following translation-equivariance variant of Proposition 4. In what follows, we say a mapping $\psi : \mathbb{M} \to \mathbb{M} \times \mathbb{M}$ is a *translation-equivariant isomorphism* from $r$ to $(r, s)$ if it satisfies all the requirements of an isomorphism, except that it may not commute with all isometries of $\mathbb{R}^d$, but only all translations of $\mathbb{R}^d$.

PROPOSITION 15. *Fix $s > 0$. There exists a finitary translation-equivariant isomorphism $\psi_s$ from $r$ to $(r, s)$, independent of $r > 0$. Furthermore, the map $\psi_s$ has finitary inverse.*

We use the following lemma in the proof of Proposition 15.

LEMMA 16. *Let $X$ be a random variable taking values in the measurable space $(A, \mathcal{A})$ and let $\Gamma(X) = \{\Gamma(X)_i\}_{i \in \mathbb{N}}$ be a random Borel partition of $\mathbb{R}^d$ which depends on $X$. Let $g : [0, 1] \times A \times \mathbb{N} \to \mathbb{M}$ be a measurable map such that if $V$ is uniformly distributed, then for all $a \in A$ and $i \in \mathbb{N}$, we have that $g(V, a, i)$ is a Poisson point process on $\Gamma(a)_i$ with intensity $s$. Let $U = \{U_i\}_{i \in \mathbb{N}}$ be a collection of independent uniform random variables independent of $X$. Then*

$$F(X, U) := \sum_{i \in \mathbb{N}} g(U_i, X, i),$$

*is a Poisson point process on $\mathbb{R}^d$ with intensity $s$ and $F(X, U)$ is independent of $X$.*

PROOF. Let $Q$ be the law of $X$ and $\Lambda$ be the law of $U$. Since $X$ is independent of $U$, for measurable $M \subseteq \mathbb{M}$ and $M' \in \mathcal{A}$, setting $L := \mathbb{P}(F(X, U) \in M, X \in M')$, we have

$$L = \int \int \mathbb{1}[F(a, u) \in M, a \in M'] \, dQ(a) \, d\Lambda(u)$$

$$= \int \int \mathbb{1}[F(a, u) \in M] \mathbb{1}[a \in M'] \, dQ(a) \, d\Lambda(u)$$

$$= \int \left[\int \mathbb{1}[F(a, u) \in M] \, d\Lambda(u)\right] \mathbb{1}[a \in M'] \, dQ(a).$$

By Remark 5 and the assumption on $g$, we have that $F(a, U)$ is a Poisson point process on $\mathbb{R}^d$ with intensity $s$ for all $a \in A$. Thus,

$$\mathbb{P}(F(X, U) \in M, X \in M') = P_s(M) \int \mathbb{1}[a \in M'] \, dQ(a)$$

(3.1)
$$= P_s(M) Q(M'),$$

which establishes the desired independence; setting $M' = A$ in (3.1) gives us that $F(X, U)$ is a Poisson point process on $\mathbb{R}^d$ with intensity $s$. $\square$

In the proof of Proposition 15, we construct a translation-equivariant isomorphism between a Poisson point process of intensity $r$ and a product of Poisson point processes of intensities $r$ and $s$. If we invert what we have done to this product, we obtain the original Poisson point process of intensity $r$. We will apply the inverse map to the permuted objects to obtain a Poisson point process of the desired intensity $s$. Since the objects are independent, this operation is well defined. This is the essential idea for the proof of Theorem 2.

PROOF OF PROPOSITION 15. Let $X$ be a Poisson point process on $\mathbb{R}^d$ of intensity $r$ with special globes $\beta_i$, each with center $c_i$ and unique point of the process $x_i$. Let $v_i$ be the cell with site $c_i$. Let $\{U_i\}_{i \in \mathbb{N}}$ be a sequence of independent uniform random variables independent of $X$. Let $\alpha$ be the map from Corollary 10 so that

$$(3.2) \qquad (X|_{\mathcal{S}^*(X)^c}, \mathcal{S}^*(X), \{\alpha(X, \beta_i)\}_{i \in \mathbb{N}}) \stackrel{d}{=} (X|_{\mathcal{S}^*(X)^c}, \mathcal{S}^*(X), \{U_i\}_{i \in \mathbb{N}}).$$

Let $\{\pi_{(A,r)}\}_{A \in \mathcal{K}}$ be the maps from Proposition 14. We write $\pi_{(v_i, s)} = \pi_{v_i}$ for simplicity, but the intensity switch is crucial to the proof.

By Proposition 14, for each $i \in \mathbb{N}$, $\pi_{v_i}(U_i)_1 \stackrel{d}{=} \pi_{v_i}(\alpha(X, \beta_i))_1$ is a Poisson point process of intensity $s$ on $v_i$, and $\pi_{v_i}(U_i)_2 \stackrel{d}{=} \pi_{v_i}(\alpha(X, \beta_i))_2$ is uniformly distributed. We need to make a slight modification to the first composition to satisfy translation-equivariance. Define the composition $\pi'_{v_i} : \mathbb{M} \times \mathcal{F} \to \mathbb{M} \times [0, 1]$ so that

$$\pi'_{v_i}(X, \beta) = \big(\pi'_{v_i}(X, \beta)_1, \pi'_{v_i}(X, \beta)_2\big)$$
$$:= \big(t_{c_i} \circ \pi_{t_{c_i}^{-1}(v_i)}(\alpha(X, \beta))_1, \pi_{v_i}(\alpha(X, \beta))_2\big),$$

where $t_{c_i}$ denotes translation by $c_i \in \mathbb{R}^d$. Still, $\pi'_{v_i}(X, \beta_i)_1$ is a Poisson point process of intensity $s$ on $v_i$. Shifting each cell to the origin, generating a Poisson point process, and shifting each cell back to center $c_i$ ensures that the generation depends on the shape of the cell but not its location. Define $\pi'(X)$ via

$$\pi'(X)|_{v_i} = \pi'_{v_i}(X, \beta_i)_1.$$

Let $b$ be defined as in Lemma 7. Recall that $b$ provides an isomorphism from a uniform random variable on $\bar{B}(\mathbf{0}, 1)$ to a uniform random variable. Let $b^{-1}$ be the inverse of $b$. Define $\mathcal{R} : \mathbb{M} \times \mathcal{F} \to \mathbb{R}^d$ so that

$$\mathcal{R}(X, \beta) = \begin{cases} t_{c_i} \circ b^{-1}\big(\pi'_{v_i}(\alpha(X, \beta_i))_2\big) & \text{if } \beta = \beta_i, \\ \mathbf{0} & \text{otherwise.} \end{cases}$$

Each $\mathcal{R}(X, \beta_i)$ is uniformly distributed in $\bar{B}(c_i, 1)$.

Set
$$X' := X - \sum_{i \in \mathbb{N}} \delta[x_i] + \sum_{i \in \mathbb{N}} \delta[\mathcal{R}(X, \beta_i)],$$
so the points in the special globes are resampled. It follows from Corollary 9, Proposition 14, and (3.2) that $X'$ is a Poisson point process of intensity $r$ on $\mathbb{R}^d$.

Define $\psi_s : \mathbb{M} \to \mathbb{M} \times \mathbb{M}$ so that
$$\psi_s(X) = (\psi_s(X)_1, \psi_s(X)_2) := (X', \pi'(X)).$$
The maps given by Proposition 14 produce a Poisson point process in the first coordinate that is independent of the uniform random variable given in the second coordinate which is used in the resampling. Define $\bar{\pi} : [0, 1] \times \mathbb{M} \times \mathbb{N} \to \mathbb{M}$ so that
$$\bar{\pi}(U_i, X, i) = \left(t_{c_i} \circ \pi_{t_{c_i}^{-1}(v_i)}(U_i)\right)_1.$$
Thus by (3.2),
$$(3.3) \qquad (X', \pi'(X)) = \left(X', \sum_{i \in \mathbb{N}} \pi'_{v_i}(X, \beta_i)_1\right) \stackrel{d}{=} \left(X', \sum_{i \in \mathbb{N}} \bar{\pi}(U_i, X, i)\right).$$

Recall that the special Voronoi cells $\mathcal{V}^*(X) = \{v_i\}_{i \in \mathbb{N}}$ give a random partition of $\mathbb{R}^d$, and by definition of $X'$ we have $\mathcal{V}^*(X) = \mathcal{V}^*(X')$, so that the random partition depends only on $X'$; in addition, the processes $X'$ and $X$ have the same special globes $\beta_i$ and centers $c_i$, so that $\bar{\pi}(U_i, X, i) = \bar{\pi}(U_i, X', i)$. Thus Lemma 16 with (3.3) gives that $\pi(X)$ is Poisson point process of intensity $s$ on $\mathbb{R}^d$ that is independent of $X'$. We already have that $X'$ is a Poisson point process of intensity $r$ on $\mathbb{R}^d$, so we have verified that $\psi_s(X)_1$ is a Poisson point process of intensity $r$ independent of $\psi_s(X)_2$, a Poisson point process of intensity $s$.

Injectivity of $\psi_s$ follows from the fact the components of $\psi_s$ are each injective; in particular, the mappings from Proposition 14 are injective and the mapping $b$ from Lemma 7 is bijective.

Next, we verify translation-equivariance. Let $\tau \in G$ be a translation. Note that the map $\alpha$ in Corollary 10 satisfies translation-invariance:
$$\alpha(\tau X, \tau \beta_i) = \alpha(X, \beta_i).$$
The map $\pi_{v_i}$ from Proposition 14 also satisfies translation-invariance in the second coordinate: $\pi_{v_i}(\cdot)_2 = \pi_{\tau v_i}(\cdot)_2$ so that
$$\pi'_{\tau v_i}(\alpha(\tau X, \tau \beta_i))_2 = \pi'_{v_i}(\alpha(X, \beta_i))_2.$$
Further, the set $\{\mathcal{R}(X, \beta_i)\}_{i \in \mathbb{N}}$ is translation-equivariant:
$$\mathcal{R}(\tau X, \tau \beta_i) = t_{\tau c_i} \circ b^{-1}\left(\pi'_{\tau v_i}(\alpha(\tau X, \tau \beta_i))_2\right)$$
$$= \tau \circ t_{c_i} \circ b^{-1}\left(\pi'_{v_i}(\alpha(X, \beta_i))_2\right) = \tau \circ R(X, \beta_i).$$

Thus in the first coordinate of $\psi_s(X)$ we have

$$\psi_s(\tau X)_1 = \tau X' = \tau X - \sum_{i \in \mathbb{N}} \tau \delta[x_i] + \sum_{i \in \mathbb{N}} \delta[\mathcal{R}(\tau X, \tau \beta_i)] = \tau \circ \psi_s(X)_1,$$

and in the second coordinate, we have

$$\begin{aligned}
\psi_s(\tau X)_2 &= \sum_{i \in \mathbb{N}} \pi'_{\tau v_i}\big(\alpha(\tau X, \tau \beta_i)\big)_1 \\
&= \sum_{i \in \mathbb{N}} t_{\tau c_i} \circ \pi_{t_{\tau c_i}^{-1}(\tau v_i)}\big(\alpha(\tau X, \tau \beta_i)\big)_1 \\
&= \sum_{i \in \mathbb{N}} \tau \circ t_{c_i} \circ \pi_{t_{c_i}^{-1}(v_i)}\big(\alpha(X, \beta_i)\big)_1 = \tau \circ \psi_s(X)_2.
\end{aligned}$$

In order to show that $\psi_s$ and its inverse are finitary, we note that if $v$ is a special Voronoi cell under $\mu$, then by construction the coordinates $\psi_s(\mu)_1$ and $\psi_s(\mu)_2$ restricted to $\mu$ are completely determined by $\mu$ restricted to $v$. Hence that both $\psi_s$ and $\psi_s^{-1}$ are finitary follows from Remark 11 and Lemma 12. $\square$

A translation-equivariant version of Theorem 2 now follows from Proposition 15, in exactly the same way Theorem 2 follows from Proposition 4. See Section 3.3 for the proof of the Theorem 2.

3.2. *Tools for isometry-equivariance.* As we move toward a proof of Proposition 4, it is helpful to recall the structure of $G$, the group of isometries of $\mathbb{R}^d$. We may write any $g \in G$ uniquely as $\tau \circ \rho$ for some translation $\tau \in G$ and orthogonal transformation (i.e., rotation or reflection) $\rho \in G$. Additionally, any translation $\tau$ corresponds to shifting by a unique point in $\mathbb{R}^d$ and any orthogonal transformation $\rho$ has a unique representation by some orthogonal matrix. We use these references interchangeably.

Let $\psi_s$ be defined as in the proof of Proposition 15. While $\psi_s$ is translation-equivariant, isometry-equivariance fails in both the first and second coordinate. Let $X$ be a Poisson point process with intensity $r$ on $\mathbb{R}^d$. Then $\psi_s(X)_1$ is the same, identical to $X$ except within the special globes, where each $x_i$ has been replaced with $\mathcal{R}(X, \beta_i)$. Recall this replacement relies on the polar coordinates of the value of $t_{c_1}^{-1}(x_i)$; hence $\psi_s(X)_1$ cannot be isometry-equivariant. In the second coordinate we have a Poisson point process $\psi_s(X)_2$ with intensity $s$ on $\mathbb{R}^d$ amalgamated from processes within each cell of $\mathcal{V}^*(X)$. We accounted for translates of cells but not orthogonal transformations in the definition of $\pi'$. Example 17 illustrates how isometry-equivariance could fail.

EXAMPLE 17 (Rotating the unit square with an interior point). Let $C_{(0,0)}$ be as in Example 6, so that $C_{(0,0)}$ is the unit square in $\mathbb{R}^2$ with lower left endpoint at the origin. Let $U_1$ and $U_2$ be uniformly distributed. Suppose we populate $C_{(0,0)}$ with a single point $Z = (U_1, U_2)$. Consider the event where $Z = (0.25, 0.75)$. Let $\rho$ be a orthogonal transformation by 90 degrees clockwise around the origin. Then $\rho(\{C_{(0,0)}, Z\})$ yields the square $C_{(0,-1)}$ with the point $\rho(Z) = (0.75, -0.25)$ inside.

However, if $C_{(0,0)}$ is first rotated by $\rho$ and then populated by the same method as in Example 6, we would obtain the point $Z + (0, -1) = (0.25, -0.25)$ inside $C_{(0,-1)}$.

To avoid any ambiguity, we make sure the output Poisson point process $\psi_s(X)_2$ on any cell does not depend on the orientation of that cell. We take a similar approach as we did for translations, but we need more machinery, also found in [6], which we now introduce. We will associate an isometry—itself equivariant under isometries—with each special globe and hence cell.

LEMMA 18. *Let $X$ be a Poisson point process on $\mathbb{R}^d$ with intensity $r$. Then the following statements hold $P_r$-almost surely.*

(a) *Distances from points of $X$ to a fixed point in $\mathbb{R}^d$ are unique.*
(b) *Inter-point distances of $X$ are unique.*
(c) *Any $d$ points of $X$ form a basis for $\mathbb{R}^d$.*

See [6], Lemma 14, for a proof. The statements essentially follow from elementary properties of Poisson point process.

Let $X$ be a Poisson point process on $\mathbb{R}^d$ of intensity $r$. Recall that the halo of a special globe $\beta_i$ contains more than $d$ points of $X$ $P_r$-almost surely. We define the *d-tag* of $\beta_i$ to be the matrix $H_i$ composed of the following $d$ columns, defined inductively. Denote the $j$th column of $H_i$ by $H_i^j$. Consider the two mutually closest points of $X$ in the halo of $\beta_i$. Of these two points, call the one closest to the center of the globe $h_i^1$. Set the first column $H_i^1 = (h_i^1 - c_i)^T$. For $1 < j \le d$ let $h_i^j$ be the point of $X$ closest to $h_i^{j-1}$ not equal to $h_i^\ell$ for any $1 \le \ell \le j-1$, and set $H_i^j = (h_i^j - c_i)^T$. Define the *d-tagging function* to be the map $\mathcal{H} : \mathbb{M} \times \mathcal{F} \to \mathbb{R}^{d \times d}$ so that

$$\mathcal{H}(X, \beta) = \begin{cases} H_i & \text{whenever } \beta = \beta_i, \\ \text{the zero matrix} & \text{otherwise.} \end{cases}$$

By its definition, each $\mathcal{H}(X, \beta_i)$ depends only on $(\mathcal{S}^*(X), X|_{\mathcal{S}^*(X)^c})$, and the set $\{\mathcal{H}(X, \beta_i)\}_{i \in \mathbb{N}}$ is equivariant under isometries. Moreover, we have that $H_i$ is nonsingular $P_r$-almost surely by Lemma 18. We use the QR factorization of $H_i$ to find our desired isometry.

Any real square matrix $A$ can be factored into a product of an orthogonal matrix $Q$ and an upper triangular matrix $R$. If $A$ is nonsingular and we require the diagonal entries of $R$ to be positive, then the factorization is unique. For a nonsingular matrix $A$ we refer to this as the *unique QR factorization of A* (see [13], Chapter 1 for details). In particular, we may write each d-tag $H_i$ as its unique QR factorization $P_r$-almost surely, which we denote by $Q_i R_i$ in the event it exists. We call $R_i$ the *upper triangular matrix* for the special globe $\beta_i$. Note $Q_i^T H_i = R_i$. We call the unique isometry that yields $R_i - Q_i^T(c_i)$ when applied to $H_i$ the *fixing isometry* for the special globe $\beta_i$. Thus we define $\sigma : \mathbb{M} \times \mathcal{F} \to G$ so that

$$\sigma(X, \beta) = \begin{cases} t^{-1}_{Q_i^T(c_i)} \circ Q_i^T & \text{whenever } \beta = \beta_i, \\ \text{the zero matrix} & \text{otherwise.} \end{cases}$$

Regarding notation, we write $Q_i^T(c_i)$ rather than $Q_i^T(c_i^T)$ for convenience, and

$$t^{-1}_{Q_i^T(c_i)} \circ Q_i^T$$

represents the translation applied to each (transposed) column of $Q_i^T$. We follow this convention throughout the section.

In the following lemma, we prove that the fixing isometry $\sigma(X, \beta_i)$ designates an isometry-invariant basis centered at the origin for the special cell $v_i$.

LEMMA 19. *Let $\mu \in \mathbb{M}$ such that the inter-point distances of $\mu$ are unique. For each special globe $\beta_i$ of $\mu$ and its d-tag $H_i$, the upper triangular matrix $R_i$ and the fixing isometry $\sigma(\mu, \beta_i)$ have the following properties.*

(a) *Each $R_i$ and $\sigma(\mu, \beta_i)$ depend only on $(\mathcal{S}^*(\mu), \mu|_{\mathcal{S}^*(\mu)^c})$.*

(b) *For each $i \in \mathbb{N}$, and for $g \in G$, the globe $\beta_i$ (under $\mu$) and the globe $g\beta_i$ (under $g\mu$) share the same upper triangular matrix; that is, the matrix $R_i$ is isometry-invariant.*

PROOF. That $R_i$ and $\sigma(\mu, \beta_i)$ depend only on $(\mathcal{S}^*(\mu), \mu|_{\mathcal{S}^*(\mu)^c})$ follows from the definitions of $R_i$ and $\sigma_i$. Let $g \in G$. By $gH_i$ we are referring to $\mathcal{H}(g\mu, g\beta_i)$, so we have that the $j$th column of $gH_i$ is $(gh_i^j - gc_i)^T$. Now, for some translation $\tau$ and orthogonal transformation $\rho$, we have $g = \tau \circ \rho$. Then

$$(gH_i)^j = (\rho(h_i^j) + \tau - \rho(c_i) - \tau)^T = (\rho(h_i^j - c_i))^T$$

and $gH_i = \rho(H_i) = \rho Q_i R_i$. The matrix $\rho Q_i$ is orthogonal, so $\rho Q_i R_i$ is the unique QR factorization of $gH_i$ and $R_i$ is the upper triangular matrix for $g\beta_i$. □

We use the fixing isometry together with the map $b$ from Lemma 7 to build an isometry-invariant version of the map $\alpha$ from Corollary 10. Since the fixing isometry depends on $(\mathcal{S}^*(X), X|_{\mathcal{S}^*(X)^c})$ we need to ensure the new output uniform random variables remain independent of $(\mathcal{S}^*(X), X|_{\mathcal{S}^*(X)^c})$. We do so with Lemma 20.

LEMMA 20. *Let $B = \bar{B}(\mathbf{0}, 1) \subseteq \mathbb{R}^d$ and let $U[B]$ be a uniform random variable on $B$. Let $\Theta$ be a random orthogonal transformation of $\mathbb{R}^d$ such that $U[B]$ and $\Theta$ are independent. Then $\Theta$ and $\Theta(U[B])$ are independent.*

PROOF. Let $S \in \mathcal{B}(\mathbb{R}^{d \times d})$ and let $A \in \mathcal{B}(B)$. Denote the law of $\Theta$ by $Q$ and the law of $U[B]$ by $L$. By the assumed independence, we have

$$\mathbb{P}(\Theta \in S, \Theta(U[B]) \in A) = \int \int \mathbb{1}[\theta \in S, \theta u \in A] \, dQ(\theta) \, dL(u)$$

$$= \int \int \mathbb{1}[\theta \in S] \mathbb{1}[u \in \theta^{-1}(A)] \, dQ(\theta) \, dL(u)$$

$$= \int \left[ \int \mathbb{1}[u \in \theta^{-1}(A)] \, dL(u) \right] \mathbb{1}[\theta \in S] \, dQ(\theta)$$

$$= \int \mathbb{P}(U[B] \in \theta^{-1} A) \mathbb{1}[\theta \in S] \, dQ(\theta).$$

Since Lebesgue measure is invariant under isometries,

$$\mathbb{P}(\Theta \in S, \Theta(U[B]) \in A) = \mathbb{P}(U[B] \in A) \int \mathbb{1}[\theta \in S] \, dQ(\theta)$$

(3.4)
$$= \mathbb{P}(\Theta \in S) \mathbb{P}(U[B] \in A).$$

Taking $S = \mathbb{R}^{d \times d}$ in (3.4), we obtain

$$\mathbb{P}(\Theta(U[B]) \in A)) = \mathbb{P}(U[B] \in A);$$

hence (3.4) also yields the required independence. □

We have assembled all the necessary pieces for the proof of Proposition 4. Indeed, much of the proof will be similar to the proof of Proposition 15. We make modifications via the fixing isometries.

PROOF OF PROPOSITION 4. Let $X$ be a Poisson point process with intensity $r$ on $\mathbb{R}^d$ with special globes $\beta_i$, each with center $c_i$ and unique point $x_i$, and special Voronoi cells $v_i$. Each special globe $\beta_i$ has d-tag $H_i$ with almost surely unique QR factorization $Q_i R_i$ and fixing isometry $\sigma(\mu, \beta_i)$. Let $b$ be defined as in Lemma 7. First, we modify the map $\alpha$ from Corollary 10. Define $\alpha' : \mathbb{M} \times \mathcal{F} \to [0, 1]$ so that

$$\alpha'(X, \beta) = \begin{cases} b(\sigma(X, \beta_i)(x_i)) & \text{if } \beta = \beta_i \text{ and } H_i \text{ is nonsingular,} \\ 0 & \text{otherwise.} \end{cases}$$

Note $\sigma(X, \beta_i)(x_i) = Q_i^T(x_i) - Q_i^T(c_i) = Q_i^T(x_i - c_i)$.

Let $\{\pi_{(A,r)}\}_{A \in \mathcal{K}}$ be the family of maps from Proposition 14. As before, it is important to note the intensity although we will write $\pi_{(v_i,s)} = \pi_{v_i}$. Define $\pi'_{v_i}$:

$\mathbb{M} \times \mathcal{F} \to \mathbb{M} \times [0, 1]$ so that

$$\pi'_{v_i}(X, \beta) = (\pi'_{v_i}(X, \beta)_1, \pi'_{v_i}(X, \beta)_2)$$
$$:= \left(\sigma(X, \beta)^{-1} \circ \pi_{\sigma(X,\beta)(v_i)}(\alpha'(X, \beta))_1, \pi_{\sigma(X,\beta)(v_i)}(\alpha'(X, \beta))_2\right),$$

where by definition

$$\sigma(X, \beta)^{-1}(z) = \left(t^{-1}_{Q_i^T(c_i)} \circ Q_i^T(z)\right)^{-1} = Q_i(z) + c_i$$

for any $z \in \mathbb{R}^d$. Set $\pi'(X)|_{v_i} := \pi'(X, \beta_i)_1$.

For the resampled points of $X$, define $\mathcal{R} : \mathbb{M} \times \mathcal{F} \to \mathbb{R}^d$ so that

$$\mathcal{R}(X, \beta) = \begin{cases} \sigma(X, \beta_i)^{-1} \circ b^{-1}\left(\pi'_{v_i}(\alpha(X, \beta_i))_2\right) & \text{if } \beta = \beta_i, \\ \mathbf{0} & \text{otherwise.} \end{cases}$$

Now set

$$X' := X - \sum_{i \in \mathbb{N}} \delta[x_i] + \sum_{i \in \mathbb{N}} \delta[\mathcal{R}(X, \beta_i)]$$

and $\psi_s(X) = (\psi_s(X)_1, \psi_s(X)_2) := (X', \pi'(X))$. We will now verify that $\psi_s$ satisfies the required properties, emphasizing the differences with the translation-equivariant case of Proposition 15.

Corollary 10 together with Lemma 20 imply that

$$(3.5) \qquad (X|_{\mathcal{S}^*(X)^c}, \mathcal{S}^*(X), \{\alpha'(X, \beta_i)\}_{i \in \mathbb{N}}) \stackrel{d}{=} (X|_{\mathcal{S}^*(X)^c}, \mathcal{S}^*(X), \{U_i\}_{i \in \mathbb{N}}).$$

Let $\bar{\pi} : [0, 1] \times \mathbb{M} \times \mathbb{N} \to \mathbb{M} \times [0, 1]$ be so that

$$\bar{\pi}(U_i, X, i) := \left(\sigma(X, \beta_i)^{-1} \circ \pi_{\sigma(X,\beta)(v_i)}(U_i)_1, \pi_{\sigma(X,\beta)(v_i)}(U_i)_2\right).$$

From (3.5) and the independence properties of the mappings from Proposition 14, we again have

$$(3.6) \qquad (X', \pi'(X)) = \left(X', \sum_{i \in \mathbb{N}} \pi'_{v_i}(X, \beta_i)_1\right) \stackrel{d}{=} \left(X', \sum_{i \in \mathbb{N}} \bar{\pi}(U_i, X, i)_1\right).$$

In addition to noting that $\mathcal{V}^*(X) = \mathcal{V}^*(X')$ and that $X$ and $X'$ have the same special globes $b_i$ and centers $c_i$, we note that by Lemma 19 the fixing isometries for $X'$ and $X$ are the same. Thus, $\bar{\pi}(U_i, X', i) = \bar{\pi}(U_i, X, i)$. Hence, Lemma 16 and (3.6) give that $\pi'(X)$ is a Poisson point process of intensity $s$ on $\mathbb{R}^d$ that is independent of $X'$.

Similarly, from (3.5) and Proposition 14 we have

$$X' \stackrel{d}{=} X - \sum_{i \in \mathbb{N}} \delta[x_i] + \sum_{i \in \mathbb{N}} \delta\left[\sigma(X, \beta_i)^{-1} \circ b^{-1}(\bar{\pi}(U_i, X, i)_2)\right].$$

Hence Corollary 9 and another application of Lemma 20 give that $X'$ is a Poisson point process of intensity $r$ on $\mathbb{R}^d$.

Next, we verify isometry-equivariance, which follows by construction. We claim $\{\alpha'(X, \beta_i)\}_{i \in \mathbb{N}}$ is isometry-invariant. Let $g \in G$. By Lemma 19, for any $z \in \beta_i$, we have $\sigma(gX, g\beta_i)(gz) = \sigma(X, \beta_i)(z)$. Thus we have $\alpha'(g\mu, g\beta) = \alpha'(\mu, \beta)$.

The set $\{\mathcal{R}(X, \beta_i)\}_{i \in \mathbb{N}}$ is isometry-equivariant: for any isometry $g = \tau \circ \rho$ and $z \in \mathbb{R}^d$ we have

$$\sigma(gX, g\beta_i)^{-1}(z) = \left(\tau^{-1}_{(\rho Q_i)^T(gc_i)} \circ (\rho Q_i)^T(z)\right)^{-1}$$
$$= \rho Q_i(z) + gc_i$$
$$= \rho(Q_i(z) + c_i) + \tau$$
$$(3.7) \qquad = g\big(\sigma(X, \beta_i)^{-1}(z)\big).$$

Then since $\{\alpha'(X, \beta_i)\}_{i \in \mathbb{N}}$ is isometry-invariant and

$$(3.8) \qquad \sigma(gX, g\beta_i)(gv_i) = \sigma(X, \beta_i)(v_i)$$

again by Lemma 19, it follows that $\{\mathcal{R}(X, \beta_i)\}_{i \in \mathbb{N}}$ is isometry-equivariant.

Isometry-equivariance in the first coordinate follows from isometry-equivariance of $\{\mathcal{R}(X, \beta_i)\}_{i \in \mathbb{N}}$. In the second coordinate, by (3.7), (3.8), and the isometry-invariance of $\{\alpha'(X, \beta_i)\}_{i \in \mathbb{N}}$ we have

$$\pi'_{gv_i}(gX, g\beta_i)_1 = \sigma(gX, g\beta_i)^{-1} \circ \pi_{\sigma(gX, g\beta_i)(gv_i)}\big(\alpha'(gX, g\beta_i)\big)_1$$
$$= g \circ \sigma(X, \beta_i)^{-1} \circ \pi_{\sigma(X, \beta_i)(v_i)}\big(\alpha'(X, \beta_i)\big)_1$$
$$= g \circ \pi'_{v_i}(X, \beta_i)_1.$$

The injectivity of $\psi_s$ follows from the fact that each of its component parts are injective. That the map $\psi_s$ and its inverse are finitary follows exactly the same argument as in Proposition 15 once we note that by Remark 11, each fixing isometry $\sigma(X, \beta_i)$ depends only on information found within the Voronoi cell $v_i$. □

3.3. *Proof of the main theorem.* Now, we have the map $\psi_s$ with all of our desired properties, except we have a surplus Poisson point process of intensity $r$. We are ready to prove the main theorem.

PROOF OF THEOREM 2. Fix $r, s > 0$. Let $X$ be a Poisson point process with intensity $r$ on $\mathbb{R}^d$ and let $\psi_s$ and $\psi_r$ be defined as in Proposition 4. Then $\psi_s(X)_1$ is a Poisson point process with intensity $r$ (which differs from $X$ only within their special globes $\beta_i$), and $\psi_s(X)_2$ is a Poisson point process with intensity $s$ independent of $\psi_s(X)_1$.

Since $\psi_s(X)_1$ and $\psi_s(X)_2$ are independent we may consider $\psi_r^{-1}$ applied to their permutation. Set $Y = \psi_r^{-1}(\psi_s(X)_2, \psi_s(X)_1)$. Then $Y$ is a Poisson point process with intensity $s$ (which differs from $\psi_s(X)_2$ only within the special globes of $\psi_s(X)_2$ and $Y$).

For $\mu \in \mathbb{M}$ set $\phi(\mu) := \psi_r^{-1}(\psi_s(\mu)_2, \psi_s(\mu)_1)$. By Proposition 4, we have that $\phi$ is a finitary isomorphism from $r$ to $s$. □

**4. Concluding remarks.** We conclude with a constructive proof of Lemma 13, a comment on the impossibility of finitary isomorphism with a fixed coding window, an application of Theorem 2 and Proposition 4 to products of Poisson systems, and a question on the property of source-universality.

4.1. *Uniform random variables on bounded and convex polytopes.* Let $v_1, \ldots, v_{d+1} \in \mathbb{R}^d$ be affinely independent, so that $v_2 - v_1, \ldots, v_{d+1} - v_1$ are linearly independent. We say $A$ is the *d-simplex* determined by $v_1, \ldots, v_{d+1}$, which we refer to as the *vertices* of $A$, if $A$ is the convex hull of its vertices. Note every point of $A$ can be uniquely expressed as a convex combination of $v_1, \ldots, v_{d+1}$.

Let $U_1, \ldots, U_d$ be uniform random variables. The *order statistics* for $U_1, \ldots, U_d$ are the random variables

$$U_{(1)} \leq U_{(2)} \leq \cdots \leq U_{(d)},$$

where $\{U_1, \ldots, U_d\} = \{U_{(1)}, \ldots, U_{(d)}\}$; see, for example, [5], Chapter 4.4, for details on the joint and marginal distributions of order statistics. Set $U_{(0)} = 0$ and $U_{(d+1)} = 1$. The *spacings* for $U_1, \ldots, U_d$ are the $d+1$ random variables $S_i = U_{(i)} - U_{(i-1)}$ for $1 \leq i \leq d+1$.

LEMMA 21. *Let $A$ be a $d$-simplex with vertices $v_1, \ldots, v_{d+1} \in \mathbb{R}^d$ and let $U_1, \ldots, U_d$ be independent uniformly distributed random variables. Let $S_1, S_2, \ldots, S_{d+1}$ be the spacings. Then the random variable*

$$\sum_{i=1}^{d+1} v_i S_i$$

*is uniformly distributed on $A$.*

For a proof of Lemma 21 see [2], Chapter XI, Theorem 2.1. We prove Lemma 13 by decomposing a bounded convex polytope into a union of $d$-simplices and applying Lemma 21. It is not difficult to generate the required spacings in an injective way from a single uniform random variable.

LEMMA 22. *There exists a measurable and injective map $h : [0, 1] \to [0, 1]^{d+1}$ such that if $U$ is uniformly distributed, then $h(U)$ has the distribution of the spacings for $d$ independent uniform random variables.*

PROOF. Let $U$ be a uniform random variable, and let $b_n$ be the binary expansion map in (2.2). Set $U_i = b_d(U)^i$ for $1 \leq i \leq d$. Then the $U_i$ are independent and each is uniformly distributed. We will construct random variables with the same distribution as the order statistics of $U_1, \ldots, U_d$, as an injective function of $U_1, \ldots, U_d$, and thus of $U$.

Let $F_1$ be the cumulative distribution function for $U_{(1)}$. Then set

$$V_1 := F_1^{-1}(U_1) \stackrel{d}{=} U_{(1)}.$$

Let $j(u, v)$ be the joint density function for $(U_{(1)}, U_{(2)})$ and $f_1$ be the density function for $U_{(1)}$. Consider the conditional distribution function given by

$$F_{2,u}(z) := \int_0^z \frac{j(u, v)}{f_1(u)} dv.$$

We set $V_2 := F_{2,V_1}^{-1}(U_2)$, so that $(V_1, V_2)$ has joint density function $j$. We similarly define $V_1, \ldots, V_d$ so that

$$(V_1, \ldots, V_d) \stackrel{d}{=} (U_{(1)}, \ldots, U_{(d)}).$$

Note that $(V_1, \ldots, V_d)$ is an injective function of $(U_1, \ldots, U_n)$.

By definition, the spacings for $V_1, \ldots, V_d$ are $S_i = V_i - V_{i-1}$. Recall we have set $V_0 = 0$ and $V_{d+1} = 1$. Clearly, $S_0, \ldots, S_{d+1}$ is a injective function of a single uniform, and has the same distribution as the spacings for $d$ independent uniform random variables. □

We remark that in the proof of Lemma 22, if we simply generated $d$ independent uniforms and then took their order statistics, we would have an $n!$ to 1 mapping.

Lemma 22 together with Lemma 21 allows us to generate a single uniform random variable on a simplex from a uniform random variable in an injective way. Since we require additional uniform random variables on the simplex, we use the following definition and lemma. Let $S^i = (S_1^i, \ldots, S_{d+1}^i)$ be independent spacings for $d$ independent uniform random variables, for $1 \leq i \leq n$. We define the *partial order statistics* $S^{(1)}, \ldots, S^{(n)}$ of $S^1, \ldots, S^n$ by ordering the spacings lexicographically so that

$$S^{(1)} \preceq S^{(2)} \preceq \cdots \preceq S^{(n)},$$

where

$$\{S^{(1)}, \ldots, S^{(n)}\} = \{S^1, \ldots, S^n\}.$$

LEMMA 23. *Let A be a $d$-simplex, and let $n \geq 1$. There exists a measurable injection $h : [0, 1] \to \mathbb{M}$ such that if $U$ is a uniform random variable, then*

$$h(U) \stackrel{d}{=} \sum_{i=1}^n \delta[U_i[A]],$$

*where $U_1[A], \ldots, U_n[A]$ are independent and uniformly distributed on $A$.*

PROOF. Let $S^1, \ldots, S^n$ be independent spacings for $d$ independent uniform random variables and let $S^{(i)} = (S_1^{(i)}, \ldots, S_{d+1}^{(i)})$ be, for $1 \leq i \leq n$ be the partial order statistics. Let $U$ be a uniform random variable. By conditioning arguments similar to those given in the proof of Lemma 22, there exists a measurable and injective map $\bar{h} : [0, 1] \to [0, 1]^{(d+1) \times n}$ so that

$$\bar{h}(U) \stackrel{d}{=} (S^{(1)}, S^{(2)}, \ldots, S^{(n)}).$$

We denote the $i$th coordinate of $\bar{h}(U)$ as $\bar{h}(U)^i$.

Let $\triangle$ be the standard simplex with vertices the unit vectors in $\mathbb{R}^{d+1}$ and let $A$ be the $d$-simplex with vertices $v_1, \ldots, v_{d+1} \in \mathbb{R}^d$. Define $\mathcal{C} : \triangle \to A$ via

$$\mathcal{C}(s_1, \ldots, s_{d+1}) = \sum_{i=1}^{d+1} v_i s_i.$$

Note $(s_1, \ldots, s_{d+1})$ is a probability vector; also note each $S^i$ takes values in $\triangle$. By Lemma 21, the random variables $\mathcal{C}(S^1), \ldots, \mathcal{C}(S^n)$ are independent and uniformly distributed on $A$. We have

$$\{\mathcal{C}(S^1), \ldots, \mathcal{C}(S^n)\} = \{\mathcal{C}(S^{(1)}), \ldots, \mathcal{C}(S^{(n)})\}$$
$$\stackrel{d}{=} \{\mathcal{C}(\bar{h}(U)^1), \ldots, \mathcal{C}(\bar{h}(U)^n)\}.$$

Thus we set

$$h(U) = \sum_{i=1}^n \delta[\mathcal{C}(\bar{h}(U)^i)]$$

so $h(U)$ has the desired law and is injective by construction. □

Recall Lemma 13 states that for a bounded convex polytope $A \in \mathcal{K}$ there exists an injective map which sends a uniform random variable to $n$ uniform random variables on $A$.

PROOF OF LEMMA 13. Let $A$ be a bounded and convex polytope. Then we may decompose $A$ into the disjoint union of $d$-simplices $A_1, \ldots, A_k$. Let $U_1[A], \ldots, U_n[A]$ be independent and uniformly distributed on $A$ for some $n \geq 1$, and set

$$X = \sum_{i=1}^n \delta[U_i[A]].$$

We will construct a measurable function $g$ so that for a uniform random variable $U$, we have $g(U) \stackrel{d}{=} X$.

Consider the set

$$\mathcal{N}_n := \{(n_1, \ldots, n_k) \in \mathbb{N}^k : 1 \leq n_1 + \cdots + n_k = n < \infty\}.$$

Let $\{\bar{n}^1, \ldots, \bar{n}^\ell\}$ be an enumeration of the set $\mathcal{N}_n$, and let

$$\bar{N} = (X(A_1), \ldots, X(A_k)).$$

Set $p_0 = 0$ and $p_j = \sum_{i=1}^{j} \mathbb{P}(\bar{N} = \bar{n}^i)$ for $1 \leq j \leq \ell$. Note that $p_\ell = 1$.

We proceed as in the proof of Proposition 14. Define $q : [0, 1) \to \mathcal{N}_n$ and $f : [0, 1) \to [0, 1)$ piecewise so that for $x \in [p_{j-1}, p_j)$,

$$q(x) = \bar{n}^j \quad \text{and} \quad f(x) = \frac{x - p_{j-1}}{p_j - p_{j-1}}.$$

Let $U$ be a uniform random variable. Then

$$f(U) \stackrel{d}{=} U \quad \text{and} \quad q(U) \stackrel{d}{=} \bar{N},$$

with $f(U)$ and $q(U)$ independent. We apply a straightforward extension of Lemma 23. For each $\bar{n} \in \mathcal{N}_n$, there exists a measurable injection $h(\bar{n}) : [0, 1] \to \mathbb{M}^k$ such that $h(\bar{n})(U)^1, \ldots, h(\bar{n})(U)^k$ are independent and for $1 \leq j \leq k$,

$$h(\bar{n})(U)^j \stackrel{d}{=} \sum_{i=1}^{n_j} \delta[U_i[A_j]],$$

where $U_1[A_1], \ldots, U_{n_j}[A_j]$ are independent and uniformly distributed on $A_j$.

Thus

$$g(U) := \sum_{j=1}^{k} h\bigl(q(U)\bigr)\bigl(f(U)\bigr)^j$$

has the same law as $X$. □

### 4.2. Remarks on fixed coding windows.

For finitary isomorphisms of Poisson systems, it is too much to ask that the size of the coding window be fixed ahead of time. In the case of one dimension, we give the following simple argument.

PROPOSITION 24. *For Poisson point processes on $\mathbb{R}$, if $0 < r < s$, then any factor from $r$ to $s$ cannot have a coding window that is a fixed deterministic constant.*

PROOF. Towards a contradiction, suppose $\phi : \mathbb{M} \to \mathbb{M}$ is a finitary factor such that if $X$ is a Poisson process of intensity $r$, then $\phi(X)$ is a Poisson process of intensity $s$, and $\phi$ has a coding window $w$ such that $P_r$-almost surely, $w \leq M$.

For each $m > 0$, let

$$E_m := \{\mu \in \mathbb{M} : \mu(B(\mathbf{0}, m)) = 0\},$$

so that $P_r(E_m) = e^{-2mr} > 0$. For $\mu \in E_M$, we must have $\phi(\mu)|_{B(\mathbf{0},1)} = \varnothing$. Since $s > r$, we choose a positive integer $\ell$ with

$$P_r(E_{M+\ell}) = e^{-2(M+\ell)r} = e^{-2Mr}e^{-2\ell r} > e^{-2\ell s} = P_s(E_\ell).$$

Since $\phi$ is translation-equivariant, $\phi(\mu)|_{B(\mathbf{0},\ell)} = \varnothing$ for all $\mu \in E_{M+\ell}$, so that

$$P_s(E_\ell) \geq P_r(E_{M+\ell}),$$

which is absurd. $\square$

We remark that we made use of the assumption that $r < s$ in our proof of Proposition 24. However, it follows from Angel, Holroyd, and Soo [1], Corollary 4, that at least in the $\mathbb{Z}$ case, there exists a translation-equivariant factor from $s = 2$ to $r = 1$ with a fixed coding window.

4.3. *Products of Poisson systems.* We apply Theorem 2 and Proposition 4 to show that Poisson systems are finitarily isomorphic to products of Poisson systems.

THEOREM 25. *Poisson systems are finitarily isomorphic to products of Poisson systems.*

PROOF. Let $r, s, s' > 0$. By Theorem 2, let $\phi_{r,s}$ be a finitary isomorphism from $r$ to $s$. Also let $\psi_{s'}$ be the map from Proposition 4.

Consider the map $\Phi = \psi_{s'} \circ \phi_{r,s}$. If $X$ is a Poisson point process on $\mathbb{R}^d$ of intensity $r$, then $\phi_{r,s}(X)$ is one of intensity $s$, and thus Proposition 4 gives that $\psi_{s'}(\phi_{r,s}(X))$ is a pair of independent Poisson processes of intensities $s$ and $s'$. All other required properties are inherited from Theorem 2 and Proposition 4. $\square$

4.4. *Source universal isomorphisms.* An interesting question raised in a paper of Harvey, Holroyd, Peres and Romik [4] asks whether there exists a source-universal finitary isomorphism of i.i.d. processes; in the context of Poisson systems, we ask the following question.

QUESTION 1. Let $s, r, r' > 0$, where $r \neq r'$. Does there exists a single measurable map $\phi : \mathbb{M} \to \mathbb{M}$ such that $\phi$ is simultaneously a finitary isomorphism from $r$ to $s$ and $r'$ to $s$?

Note the homomorphism in Theorem 3 is finitary and source-universal, and the map in the proof of Proposition 4 is also finitary and source-universal. However, the isomorphism in Theorem 2 is not, as constructed, source-universal.

**Acknowledgments.** We thank the referee for helpful comments and Benjy Weiss for helpful discussions.

DEPARTMENT OF MATHEMATICS
UNIVERSITY OF KANSAS
405 SNOW HALL
1460 JAYHAWK BLVD.
LAWRENCE, KANSAS 66045-7594
USA
E-MAIL: math@terrysoo.com
awilkens@ku.edu